\documentclass[11pt]{article}
\usepackage{graphicx}
\usepackage[update,prepend]{epstopdf}
\usepackage{amssymb}
\usepackage{amsmath}
\usepackage{amsfonts}
\usepackage{amsopn}
\usepackage{amsthm}

\def\ds{\displaystyle}

\def\L{{\bf L}}

\def\bfe{{\bf e}}
\def\ve{\varepsilon}

\def\P{{\cal P}}
\def\caL{{\mathcal L}}

\def\R{I\!\!R}
\def\N{I\!\!N}

\def\implies{\Longrightarrow}

\def\c{\centerline}
\def\P{{\cal P}}

\def\vsk{\vskip 4em}
\def\v{\vskip 1em}

\def\C{{\cal C}}
\def\ov{\overline}
\def\Tilde{\widetilde}

\def\bega{\begin{array}}
\def\enda{\end{array}}
\def\begi{\begin{itemize}}
\def\endi{\end{itemize}}

\def\bel{\begin{equation}\label}
\def\eeq{\end{equation}}
\def\sqr#1#2{\vbox{\hrule height .#2pt
\hbox{\vrule width .#2pt height #1pt \kern #1pt
\vrule width .#2pt}\hrule height .#2pt }}
\def\square{\sqr74}
\def\endproof{\hphantom{MM}\hfill\llap{$\square$}\goodbreak}

\renewcommand{\div}{\mbox{div}\,}

\begin{document}
\title{\bf Non-existence and Non-uniqueness for Multidimensional
 Sticky Particle Systems}
\vsk

\author{Alberto Bressan$^{(*)}$ and Truyen Nguyen$^{(**)}$\\  \,  \\
(*) Department of Mathematics, Penn State University.
University Park, PA~16802, USA.\\
(**) Department  of Mathematics, University of Akron.
Akron, OH 44325-4002, USA.
\\
e-mails:~ bressan@math.psu.edu~,~tnguyen@uakron.edu}

\maketitle
\begin{abstract}
The paper is concerned with sticky weak 
solutions to the equations of pressureless gases in two or more space dimensions.
Various initial data are constructed, showing that 
the Cauchy problem can  have (i) two distinct
sticky solutions, or (ii) no sticky solution, not even locally in time.
In both cases the initial density is smooth with compact support, while 
the initial velocity field is continuous.
\end{abstract}

\section{Introduction}
We consider the initial value problem for the equations of pressureless gases 
in several space dimensions:
\begin{equation}\label{SP}
\left\{\begin{array}{rl}
\!\!\partial_t\rho+\div( \rho v)&=0, \cr
\!\!\partial_t(\rho v)+\div\big( \rho v\otimes v\big)&=0,\enda
\right.  \qquad\quad  t\in  ]0,T[\,,\quad x\in\R^n,\eeq
\bel{ic}\rho(0,x)~=~\bar \rho(x), \qquad v(0,x)~=~\bar v(x).\eeq
The  system \eqref{SP} 
was first studied by Zeldovich \cite{Zeld} in the one-dimensional case to model the evolution of a sticky particle system. An example of a  measure-valued solution
is provided by a finite collection of particles moving with constant speed 
 in the
absence of forces. Whenever two or more particles collide,
they stick to each other as a single compound particle. The mass of the new particle
is equal to the total mass of the particles involved in the collision, while 
its velocity is determined by the conservation of momentum. 
The sticky particle system has been investigated extensively by many authors and  
is well understood in dimension  $n=1$, 
see \cite{BoJ,BGSW,BG, ERS,HW,NS,NT1,NT2}. 
In this  case it is known that, for any initial data $(\bar \rho, \bar v)$
with bounded total mass and energy,
 the Cauchy problem  \eqref{SP}-(\ref{ic})  
has a unique global entropy-admissible  weak solution $(\rho, v)$ 
(see \cite{NS} and \cite[Theorem~1.3]{NT2}). 

The present paper is concerned with the initial value problem associated with \eqref{SP} 
in space  dimension $n\geq 2$. We are interested in weak solutions to \eqref{SP} obeying the sticky particle or adhesion dynamics principle, which are the most relevant from a 
physical point of view. 
For initial data containing finitely many particles, it is easy to see that
a unique global solution exists, but it does not depend continuously on the initial data.
In the case of countably many particles, we show that both uniqueness and existence
can fail.   Indeed, we construct a Cauchy problem having exactly two solutions,
and a second Cauchy problem where no solution exists, not even locally in time.
Both these examples can be adapted to the case of $\L^\infty$  initial data.

The remainder of the paper is organized as follows.   In Section~2 we give precise
definitions of ``weak solution" and ``sticky solution" for initial data containing
countably many point masses and also for continuous mass distributions, following 
\cite{S}.
  Section~3 contains an example showing the non-uniqueness
of sticky solutions.
In  Section~4 we describe a Cauchy problem without any local solution
and explain how this counterexample relates to the (erroneous) proof of global
existence of sticky solutions proposed in \cite{S}.
Finally, in Section~5 we extend the analysis to  initial data having smooth density
and continuous velocity.
Even in this case we show that local existence and uniqueness do not hold, in general.

\section{Dynamics of sticky particles}
\label{sec:2}
\setcounter{equation}{0}
We consider a system containing countably many sticky particles, moving in 
$n$-dimensional space.
Let
$$\left.\bega{cl}x_i(t)&=~\hbox{position}\cr
v_i(t)&=~\hbox{velocity}\cr
m_i&=~\hbox{mass}\cr\enda\right\}
~~\hbox{of the $i$-th particle at time $t$.}$$ 
In a Lagrangian formulation, the state of the system is
described by countably many ODEs for the variables $x_i$. Let 
\bel{idata}x_i(0)~=~\bar x_i\,,\qquad\qquad \dot x_i(0+)~=~\bar v_i\eeq
be the initial position and the initial velocity of the $i$-th particle.
It is natural to assume that, when particles are at a same location, 
they stick together traveling with a common  speed
determined  by the conservation of momentum.   At any time $t\geq 0$, the speed of the
$i$-th particle should thus be
\bel{Vi}
\dot x_i(t)~=~V_i(t)~\doteq~\frac{\sum_{j\in J_i(t)}   m_j \bar v_j}
{\sum_{j\in J_i(t)}   m_j}\,,\eeq
where
\bel{Ji} J_i(t)~\doteq~\Big\{j\geq 1:\quad x_j(t)=x_i(t)\Big\}.\eeq
Notice that the right hand side of (\ref{Vi}) 
is well defined provided that the total mass $M=\sum_i m_i$
and  the initial energy $E\doteq \frac{1}{2}\sum_i m_i |\bar v_i|^2$
are finite.

\begin{figure}
\centering
\includegraphics[scale=0.45]{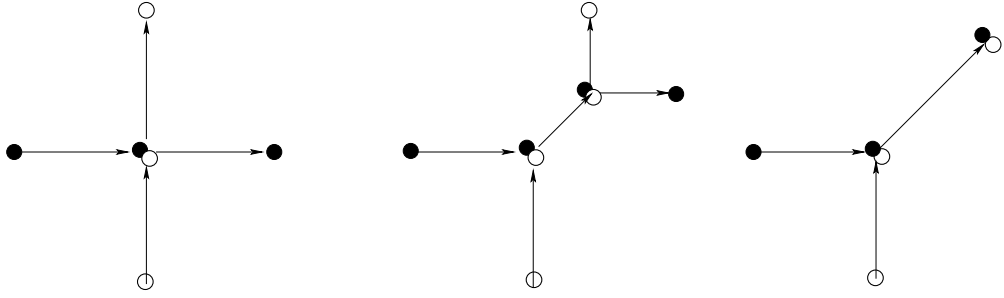}
\caption{\small  Left: a weak solution which is energy admissible but not sticky.
Center: a weak solution which is neither sticky nor energy admissible.
Right: a sticky solution.}
\label{f:z104}
\end{figure}

{\bf Definition 1.} A family of  continuous maps
$t\mapsto x_i(t)$ is a {\it weak solution} of the  equations (\ref{Vi})
with initial data (\ref{idata}) if, for every $i\geq 1$ and $t\geq 0$, one has
\bel{ws}
x_i(t)~=~\bar x_i + \int_0^t V_i(s)\, ds\,.\eeq
In addition, we say that the solution is {\it energy admissible}
if the corresponding energy
\bel{E} E(t)~\doteq~\frac12\sum_{i}m_i|\dot x_i(t+)|^2\eeq
is a bounded, non-increasing function of time.
\v
In general, even in the case of two particles, the energy-admissible weak solution
need not be unique. To achieve uniqueness (at least for finitely many particles)
one more condition must be imposed.
\v
{\bf Definition 2.}
We say that a weak solution $\{ x_i(\cdot)\,;~~i\geq 1\}$ is a {\it sticky solution}
if it satisfies the additional property
\begi
\item[{\bf (SP)}]~~ {\it
If $x_i(t_0)=x_j(t_0)$ at some time $t_0\geq 0$, then
$x_i(t)=x_j(t)$ for all $t>t_0$.}
\endi

\v
{\bf Example 1.} Given any initial data (\ref{idata}), 
the family of functions
\bel{triv}x_i(t)~=~\bar x_i + t\bar v_i\eeq
always provides a weak, energy admissible solution.    Indeed, for any $i$, the set
of times 
$$\{ t>0: \quad x_j(t)= x_i(t) ~\hbox{for some}~j\not= i\}$$
is at most countable, therefore it has measure zero.
This implies $J_i(t)=\{i\}$ for a.e.~$t\geq 0$, hence (\ref{ws}) trivially holds.
\v
{\bf Example 2.}   Consider two particles moving in the plane, with
masses $m_1=m_2=1$ and with initial
positions and velocities given by
\bel{id2}\bar x_1 = \begin{pmatrix}1\cr 0\end{pmatrix},\quad \bar x_2 = 
\begin{pmatrix}0\cr 1
\end{pmatrix},    \qquad\qquad \bar v_1 = \begin{pmatrix}0\cr 1
\end{pmatrix}, \quad  \bar v_2 = \begin{pmatrix}1\cr 0\end{pmatrix}.\eeq
As in Example 1, the maps
\bel{ex2}x_1(t)=\bar x_1+ t\bar v_1\,,\qquad\qquad x_2(t)=\bar x_2+ t\bar v_2\eeq
provide an energy-admissible weak solution, which however does not satisfies the 
stickiness assumption {\bf (SP)}.
The unique solution that satisfies {\bf (SP)} is given by \eqref{ex2} for $t\in [0,1]$, while
\bel{ex22}x_1(t)~=~x_2(t)~=~\frac{t+1}{2}\begin{pmatrix}1\cr 1\end{pmatrix}\qquad\qquad 
\hbox{for}~~~t\geq 1\,.\eeq
We remark that the weak solution \eqref{ex2} depends continuously
on the initial data, but the sticky solution does not.  Indeed, if we 
slightly perturb the initial data, say by taking 
$\bar x_1 = \begin{pmatrix}1\cr \ve\end{pmatrix}$ with $\ve\not= 0$, 
then the two particles do not collide
and (\ref{ex2}) provides the unique sticky solution to the Cauchy problem.

We also  observe that the above Cauchy problem has infinitely many 
weak solutions which are not energy admissible (and not sticky).  
Indeed, for any given time $T\geq 1$, a weak solution is defined by (\ref{ex2})
for $t\in [0,1]$, by (\ref{ex22}) for $t\in [1,T]$, and by 
\bel{ex23}
x_1(t)~=~\frac{T+1}{2}\begin{pmatrix}1\cr 1\end{pmatrix}+ 
\begin{pmatrix}0\cr t-T\end{pmatrix},\qquad
x_2(t)~=~\frac{T+1}{2}\begin{pmatrix}1\cr 1\end{pmatrix}+ 
\begin{pmatrix}t-T\cr 0\end{pmatrix}\qquad
\quad 
\hbox{for}~~~t\geq T\,.\eeq
\v
{\bf Remark 1.}  When only finitely many particles are present, is it
an easy matter to prove the global existence and  uniqueness of a sticky solution
to the Cauchy problem (\ref{idata}).  The proof  can be achieved by induction
on the number $N$ of particles.  When $N=1$ the result is trivial.
Next, assume that the result is true whenever the initial number  of particles is $<N$.    Consider an initial data consisting of exactly $N$ particles. Let $x_i(t) ~\doteq~ \bar x_i + t \bar v_i$ and 
define the first interaction time  
$$\tau~\doteq~\inf\Big\{t> 0\,;~~x_i(t)=x_j(t)\quad\hbox{for some} ~i\not= j\Big\}.$$
If $\tau=+\infty$, then 
$$x_i(t) ~=~ \bar x_i + t\bar v_i\qquad\qquad i=1,\ldots,N,$$ 
describes the unique  sticky solution. If $\tau<\infty$, then  
at time $\tau$ two or more colliding particles are lumped together
in a single compound particle with speed determined by the conservation of momentum.
This yields a new Cauchy problem, where the initial data at $t=\tau$ 
contains a number of particles
strictly less than $N$.   The result follows by induction.
\v

For initial data where the density can be an arbitrary measure, 
a general notion of sticky weak solutions for the system  \eqref{SP} 
 was introduced by Sever in \cite{S}. 
This definition is
reviewed here, and will be later used in Section~5. 
In the following we assume that the initial density and the initial 
velocity of the pressureless gas satisfy
\begin{equation}\label{initial-data}
\bar \rho\in \P_2(\R^n),\qquad\qquad\quad \bar v\in \L^2(\bar \rho).
\end{equation}
Here $\P_2(\R^n)$ is the set of all probability measures $\rho$ such that
$$\int |x|^2\, d\rho(x)~<~\infty.$$
Let $D$ be an open set in $\R^n$ with $\caL_n(D)= 1$, where $\caL_n$ is the Lebesgue measure on $\R^n$.  Regarding $y\in D$ 
as a Lagrangian coordinate, a flow will be described by a mapping $X: 
[0,T]\times D \to \R^n$ with forward time derivatives $X_t$ satisfying
\[
\sup_{0\leq t\leq T}~\int_{D}{|X_t(t,y)|^2\, dy}<\infty. 
\]
In addition to the usual Lebesgue spaces  $\L^2(D)$ and  $\L^2\big([0,T]\times D)$,
 we shall also need the following function spaces:
\begin{itemize}
\item $J(X)$ is the completion of 
$C_c^1([0,T]\times \R^n)$ with respect to the norm
\[
\|\theta\|_J~\doteq~  \Big[\int_0^T\int_{D} \theta(t, X(t,y))^2 \, dy dt\Big]^{1/2}.
\]
\item For each $t\in [0,T]$, the $t$-section  $J(X, t)$ is the completion of 
$C_c^1(\R^n)$ with respect to the norm
\[
\|\phi\|_{J,t}~ \doteq ~ \Big[\int_{D} \phi(X(t,y))^2 \, dy\Big]^{1/2}.
\]
\item $K(X)~\doteq~ \big\{\theta \circ X\,;~~ \theta \in J(X)\big\}~
\subset~ \L^2([0,T]\times D)$ ~~ 
and ~~$K(X, t)~\doteq~
 \big\{\phi \circ X(t,\cdot)\,;~~\phi \in J(X,t)\big\}~\subset~ \L^2( D)$.
\end{itemize}

It was shown in \cite[Lemma~3.1]{S} that for every $g\in \L^2( D)$ 
there exists a unique $v\in J(X)$ such that, for every $t\in [0,T]$, 
\begin{equation}\label{pre1}
v(t,\cdot)\in J(X,t),\quad  \quad\mbox{and}\quad \|v(t,\cdot)\|_{J,t}\leq \|g\|_{\L^2}.
\end{equation}
In addition, the following orthogonality relation holds:
\begin{equation}\label{pre2}
v\circ X(t,\cdot) -g \, \, \perp \, K(X,t).\eeq
Let $W(X) : L^2( D) \mapsto J(X)$ be the linear mapping 
defined by $W(X) g \doteq v$, where $v$ is the  unique  $v\in J(X)$ satisfying 
\eqref{pre1}-\eqref{pre2}. We can now recall the definition of sticky weak
solution introduced in \cite{S}.
\v
 {\bf Definition 3.} 
 A flow  mapping $X: 
[0,T]\times D \to \R^n$ provides a {\it sticky weak solution} to the Cauchy 
problem \eqref{SP}--\eqref{ic} if it satisfies the following properties.

{\bf 1 - weak solution:~} $\, X_t = (W(X) X_t(0,\cdot))\circ X \,$ in $\, 
\L^2([0,T]\times D)$.

{\bf 2 - initial data:~} For every $\phi\in J(X,0)$ one has
\begin{align*}
&\int_{D}{\phi(X(0,y))\, dy}
= \int_{\R^n}{\phi(x)\, d\bar \rho(x)},\\
&\int_{D}{ X_t(y,0) \phi(X(0,y))\, dy}
= \int_{\R^n}{\bar v(x)\phi(x)\, d\bar \rho(x)}.
\end{align*}

{\bf 3 - sticky property:~} For every $y_1,\, y_2,\, t_0$ 
such that $X(t_0,y_1)=X(t_0,y_2)$, one has
\[
X(t,y_1)=X(t,y_2)\quad \mbox{for all}\quad t\geq t_0.
\]

\v

If $X$ satisfies {\bf 1 - 2} in Definition~3, then 
\begin{equation}\label{Lagrange-Eulerian}v ~\doteq~ W(X) X_t(0,\cdot)
\qquad\hbox{and}\qquad
\rho(t,\cdot) ~\doteq~ X(t,\cdot)_{\#}\caL_n 
\end{equation}
(i.e., the push-forward of the Lebesgue measure $\caL_n$ on $D$ by the map $y\mapsto X(t,y)$)
provide a distributional solution
to  the system \eqref{SP} of pressureless gases
(see \cite[Theorem~3.2]{S} and the subsequent remark). 
Moreover, $v\circ X = X_t$ in  $\L^2([0,T]\times D)$.

\section{A Cauchy problem with two solutions}
\label{sec:3}
\setcounter{equation}{0}
We construct here
an initial configuration containing countably many 
particles in the plane, such that the Cauchy problem has two distinct 
sticky solutions.

\begin{figure}[htbp]
\centering
 \includegraphics[scale=0.45]{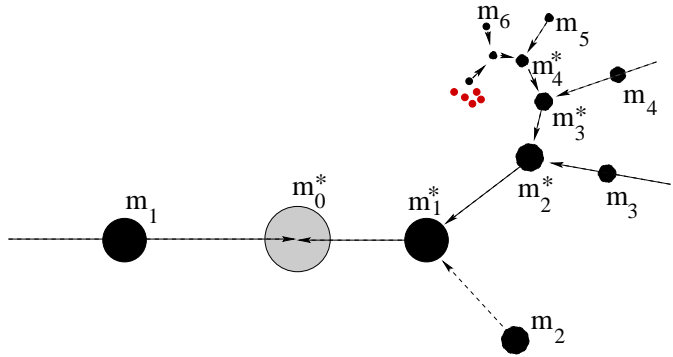}
    \caption{\small Countably many particles with masses $m_i=2^{-i}$
    collide, two at a time,
    until a single compound particle is formed.}
\label{f:z103}
\end{figure}

{\bf Example 3.} 
A solution of (\ref{ws}) will be constructed by induction, starting from the final configuration and 
going backward in time (see Figure~\ref{f:z103}).  For notational simplicity,
we denote with the same symbol a particle and its mass.  We also define the
times $t_i\doteq 2^{-i}$, $i=1,2,\ldots$
The solution is constructed by induction on the intervals $[t_{i+1},\, t_i]$.
\begi
\item For $t\geq t_1=1/2$ there is one single particle of mass $m_0^*=1$, located at the origin.

\item For   $t\in [t_2, t_1[\,$   there are two particles with  masses $m_1=m_1^*=1/2$, 
moving toward each other with opposite speed
and colliding at the origin at time $t=t_1$.

\item For   $t\in [t_3, t_2[\,$, the particle $m_1^*$ is replaced by  two equal 
particles with masses $m_2=m_2^* = 1/4$, colliding at time $t_2$.

 \c{$\ldots$}
 
\item In general,  for  $t\in [t_{i+1}, \,t_i[\,$, the particle $m_{i-1}^*$ 
is replaced by  two equal 
particles, with masses $m_i=m_i^* = 2^{-i}$, colliding at time $t_i$.
\endi

Let $ v_i$ be the   constant speed of the $i$-th particle $m_i$
   during the time interval $[0, t_i]$, i.e. before interaction. 
  Moreover, let $v_i^*$ be the speed of the lumped particle $m_i^*$ during the
  time interval $[t_{i+1}, \, t_i]$. 
Let $x_i^*\in \R^2$ be the point in the plane where the particles 
$m_i, m^*_i$ interact (at time $t_i$).
Conservation of momentum
requires
\bel{momc}v^*_{i-1} ~=~\frac{v_i + v^*_i}{2}\,.\eeq
Apart from (\ref{momc}), the speeds $v_i, v_i^*$ can be freely chosen.
We take advantage of this fact, choosing $v_i, v_i^*$ so that the following
non-intersection property holds.
\v
{\bf (NIP)} For each $t\geq 0$, the two points
$$x_i(t)~\doteq~ 
x_i^*+(t_i-t)v_i,\qquad x_i^*(t)~\doteq~ 
x_i^*+(t_i-t)v_i^*$$
do not coincide with any of the finitely many points
$$ x_j(t)~\doteq~ 
x_j^*+(t_j-t)v_j,\qquad x_j^*(t)~\doteq~ 
x_j^*+(t_j-t)v_j^*\qquad\qquad 1\leq j <i\,.$$
\v
By induction on $i=1,2,\ldots$ we thus obtain a sticky solution to the 
 particle equation (\ref{Vi})
defined for all $t\geq 0$.  The initial data consists of countably many 
particles with masses $m_i = 2^{-i}$, $i\geq 1$. 
As time progresses, these particles collide and stick to each other.  In particular 
during each time interval $[t_i, t_{i-1}[$ only $i$ distinct 
compound particles are present. 
 
We now observe that, because
of (NIP), the trivial solution \eqref{triv} with $\bar v_i ~ \doteq~ v_i$ is a sticky solution as well.
This provides a counterexample to uniqueness, for 
an initial configuration containing  countably many particles.

\v

\section{A Cauchy problem with no solution}
\label{sec:4}
\setcounter{equation}{0}
In this section we shall construct 
an initial configuration containing countably many particles
 in the plane, such that the corresponding Cauchy problem has no sticky solution,
 even locally in time.

 As a first step,
consider 
a one-dimensional configuration consisting of  countably many 
particles $x_k$, $k\geq 1$ moving along the $x$  axis.
The masses $m_k$ of these particles, and their 
initial positions $\bar x_k$  and velocities $\bar v_k$ are chosen to be 
\bel{id4}
m_k~\doteq~\alpha^k,\qquad 
\bar x_k~=~\beta^k,
\qquad \bar v_k ~=~ 1-\gamma^k.
\eeq
with $0<\alpha,\beta,\gamma <1$.   For notational convenience, we denote by
$$x_j(t)~=~ \bar x_j+ t \bar v_j\qquad\qquad j\geq 1$$
the positions of the free particles. 
Notice that the collection of particles 
$\{x_j(t)\,;~~j\geq k\}$ with masses $m_j$ has barycenter located at
\bel{bar} \bega{rl}b_k(t)&\ds\doteq~{\sum_{j\geq k} m_j(\bar x_j+ t \bar v_j)\over
\sum_{j\geq k} m_j}~=~{\sum_{j\geq k} \alpha^j[\beta^j+ t (1-\gamma^j)]\over
\sum_{j\geq k} \alpha^j}\cr\cr
&\ds=~{1-\alpha\over\alpha^k}\cdot\left[ {\alpha^k\beta^k\over1-\alpha\beta}+
{t\alpha^k\over 1-\alpha} -{t\alpha^k\gamma^k\over 1-\alpha\gamma}\right].
\enda\eeq
Call $t_{k-1}$ the time when this barycenter hits the particle $x_{k-1}$.
Solving the equation $b_k(t) = x_{k-1}(t)$ we obtain
\begin{align*}
&\beta^{k-1}-t\gamma^{k-1}~=~(1-\alpha)\left[{\beta^k\over 1-\alpha\beta}-{t\gamma^k\over 1-\alpha\gamma}\right],\\
&{1-\beta\over 1-\alpha\beta} \beta^{k-1}~=~{1-\gamma\over 1-\alpha\gamma}
\cdot t\gamma^{k-1}.
\end{align*}
Therefore,
\bel{tk-1}
t_{k-1}~=~{1-\beta\over 1-\alpha\beta}\cdot {1-\alpha\gamma\over 1-\gamma}
\cdot \left({\beta\over\gamma}\right)^{k-1}.\eeq
In particular, if we choose  $0<\beta<\gamma<1$, the sequence
of times $t_k$ will be strictly decreasing to zero.

 \begin{figure}
\centering
\includegraphics[scale=0.55]{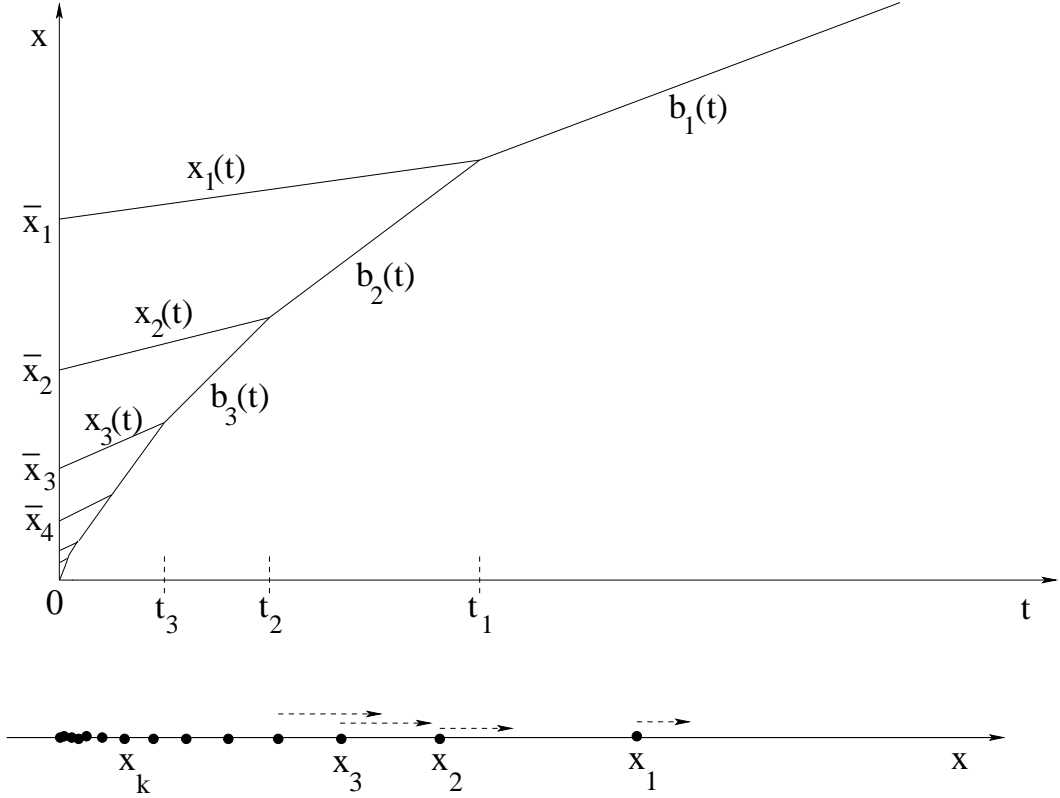}
\caption{\small  A sticky solution containing  countably many particles, 
moving on the 
$x$-axis. }
\label{f:z105}
\end{figure}

The unique solution to the one-dimensional Cauchy problem can be explicitly described 
as follows (Figure~\ref{f:z105}).
For $t\in [t_k, \, t_{k-1}[$ there are $k-1$ particles with masses $m_{k-1},\ldots, m_1$, located
at $x_{k-1}(t)<x_{k-2}(t)<\cdots<x_1(t)$, and one compound particle with 
mass $m_k^*= \sum_{j\ge k} m_j$, located at $b_k(t)$.

For future use, two lemmas will be needed.
\v
{\bf Lemma 1.} {\it If  $0<\beta<\gamma<1$
and $0<\alpha< \frac{1}{1+\beta+\gamma}$, then for every $k>1$ one has
\bel{ai}
x_{k+1}(t_{k-1})~>~x_{k-1}(t_{k-1})~=~b_k(t_{k-1})\eeq
}

{\bf Proof.}
The inequality (\ref{ai}) holds provided that
$$\beta^{k+1}-\gamma^{k+1}t_{k-1}~>~\beta^{k-1}-\gamma^{k-1}t_{k-1}.$$
An explicit computation yields
$$t_{k-1}~>~{1-\beta^2\over 1-\gamma^2}\left({\beta\over\gamma}\right)^{k-1}.$$
By \eqref{tk-1}, this is equivalent to
\begin{equation*}
{1+\beta \over 1+\gamma }~<~
 {1-\alpha\gamma\over 1-\alpha \beta}.
 \end{equation*}
Therefore, if  $0<\beta<\gamma<1$, the above inequality holds as soon as $0<\alpha< \frac{1}{1+\beta+\gamma}$. 
\endproof

{\bf Lemma 2.} {\it Let $0<\beta<\gamma<1$ and $0<\alpha< \frac{1}{1+\beta+\gamma}$ be as in Lemma~1.
Then for every $k> 1$, there exists a time $\tau_k~\in~]t_k, \, t_{k-1}[\,$
such that the following holds.

Consider any subset
$S_k'\subset S_k\doteq\{j\,;~~j\geq  k\}$, with $k\in S_k'\not= S_k$.
Then the barycenter $b_k'(\tau_k)$ of the set $\{x_j(\tau_k)\,;~~j\in S_k'\}$ satisfies
\bel{bar'}
b_k'(\tau_k)~\doteq~{\sum_{j\in S_k'} m_j x_j(\tau_k)\over
\sum_{j\in S_k'} m_j}~<~{\sum_{j\in S_k} m_j x_j(\tau_k)\over
\sum_{j\in S_k} m_j}~=~b_k(\tau_k).
\eeq
}

{\bf Proof.} Let $k> 1$ be given.  Using Lemma 1, by continuity we can find
$\tau_k~\in~]t_k, \, t_{k-1}[\,$ such that 
\[ x_{k+1}(\tau_k)~>~b_k(\tau_k).
\]
We claim that with this choice the inequalities \eqref{bar'} hold as well.
Indeed,
$$x_k(\tau_k)~<~b_k(\tau_k)~<~x_{k+1}(\tau_k)~<~x_{k+2}(\tau_k)~<~\cdots$$
Therefore, defining the set $S_k''\doteq S_k\setminus S_k'$, the corresponding
 barycenter satisfies
\bel{aai} b_k''(\tau_k)~\doteq~{\sum_{j\in S_k''} m_j x_j(\tau_k)\over
\sum_{j\in S_k''} m_j}~\geq~\min_{j\in S_k''} x_j(\tau_k)~>~b_k(\tau_k).
\eeq
Observe that
\[
b_k(\tau_k) = \theta b_k'(\tau_k) + (1-\theta) b_k''(\tau_k)
\]
for $\theta ~\doteq~ \frac{\sum_{j\in S_k'}{m_j}}{\sum_{j\in S_k}{m_j}}$.
But since $0<\theta<1$ due to $S_k'\neq \emptyset$ and $S_k' \subsetneqq S_k$, it follows from \eqref{aai} that
 $b_k(\tau_k)>b_k'(\tau_k)$.
\endproof

After these preliminaries we can describe our main counterexample.

{\bf Example 4.}  
On the plane 
$\R^2$ we shall use the canonical basis $\bfe_1= \begin{pmatrix}1\cr 0\end{pmatrix}$,
$\bfe_2= \begin{pmatrix}0\cr 1\end{pmatrix}$.
The initial configuration consists of two countable sets of particles (Figure~\ref{f:z106}).
\begi
\item A sequence of black particles moving horizontally along the $x_1$ axis.
As in (\ref{id4}), 
their masses, initial positions, and initial velocities are defined as
\bel{id6}
m_k~\doteq~\alpha^k,\qquad 
\bar x_k~=~\beta^k \bfe_1,
\qquad \bar v_k ~=~ (1-\gamma^k)\bfe_1\,.
\eeq

\item A sequence of white particles, moving vertically.   Their 
masses, initial positions, and initial velocities are chosen as
\bel{id7}
M_k~\doteq~\alpha^k,\qquad 
\ov X_k~=~b_k(\tau_k) \bfe_1 +\tau_k \bfe_2\,,
\qquad \ov V_k ~=~ -\bfe_2\,.
\eeq
\endi
We think of the white particles as bullets, knocking the black particles 
away from the $x_1$ axis.

We claim that, with this initial configuration, no sticky  solution exists.
Indeed, assume that a solution exists, and 
call $S\subset\N$ be the set of all white particles that hit a 
target, i.e.~that collide with a lumped black particle while crossing the $x_1$ axis.
Two cases can be considered, each leading to a contradiction.

CASE 1: $k\in S$ for some $k\geq 1$.   We claim that this is possible 
only if $j\notin S$ for all $j>k$.     Indeed, let $S'_k\subseteq S_k=\{j\;~j\geq k\}$ 
denote  the set of black particles which are NOT hit by some white particle
before time $\tau_k$.    If $S_k'\not= S_k$, then by Lemma~2 the barycenter
of these particles satisfies
$b'_k(\tau_k)<b_k(\tau_k)$.   Hence 
$$b'_k(\tau_k)\bfe_1~\not=~b_k(\tau_k)\bfe_1~=~X_k(\tau_k)$$
and the $k$-th white particle will not hit its target.   

On the other hand, if $j\notin S$ for all $j>k$, then none
of the black particles $x_j$ with $j\geq k+1$ is hit by white bullets.
At time $\tau_{k+1}$ the barycenter of this set of black particles
is located 
at $b_{k+1}(\tau_{k+1}) \bfe_1~=~X_{k+1}(\tau_{k+1})$.   As a result, the white particle 
$X_{k+1}$ hits its target.
We have thus proved the two implications
\bel{imp}\bega{rl}   k\in S\quad &\implies\quad j\notin S~~\hbox{for all}~~j>k,\cr\cr
j\notin S~~\hbox{for all}~~j>k\quad &\implies\quad k+1\in S,\enda\eeq
leading to a contradiction.

CASE 2: The remaining possibility is that
 $S=\emptyset$. But in this case the second implication in (\ref{imp})
 immediately yields a contradiction.
 
 {}From the above arguments, it is clear that a sticky solution does not exist,
 even locally in time.

\begin{figure}
\centering
\includegraphics[scale=0.5]{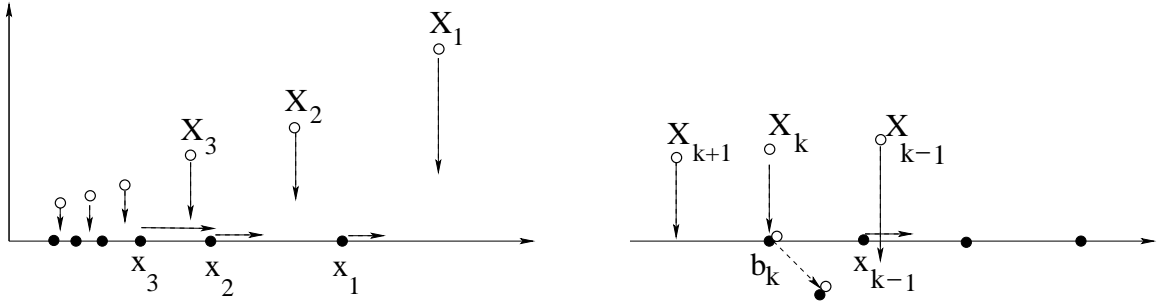}
\caption{\small Left: The initial configuration of black particles $x_k$ 
moving horizontally
and white particles $X_k$ moving vertically. Right: if the white 
particle $X_k$ scores a hit, then no other white particle $X_j$ with 
$j\not=k$
can collide with a black particle along the $x_1$ axis. }
\label{f:z106}
\end{figure}

{\bf Remark 2.} 
 Following \cite{S}, one can construct a family of weak
solutions as follows.    At each time where an interaction
occurs, two particles can either
stick together, or continue their separate motion without any change in the velocities.
The choice (sticking to each other or not) is made in order to minimize the integral
$$J_\ve~\doteq~\int_0^\infty e^{-t/\ve} E(t)\, dt$$
where $E(t)$ is the total energy at time $t$.   
It is clear that, as two interacting particles stick together, the momentum is conserved but the energy decreases.    Letting $\ve\to 0$, in \cite{S} it was claimed (but not proved) that any limit of a sequence of weak solutions which minimize $J_\ve$
should yield a solution to the sticky particle equations.
This is true in the case of finitely many particles, as suggested by intuition, 
but false in general.    In our specific example, for a given $\ve>0$ a
solution which minimizes $J_\ve$ can be described as follows.
The black particles traveling along the $x_1$ axis always stick to each other
after collision.    On the other hand, there is an integer $N = N(\ve)$ 
such that:
\begi
\item[(i)] 
For $k>N$, at time $\tau_k$ the white particle $X_k$ and  the corresponding lumped 
black particle hit each other  at $b_k (\tau_k)\bfe_1$, 
without changing their speed (i.e., without sticking).
\item[(ii)] At time $\tau_N$, the white particle $X_N$ hits  the corresponding black particle at
$b_N (\tau_N)\bfe_1$ and sticks to it, knocking it away from the $x_1$ axis.
\item[(iii)] For $k<N$, the remaining white particles do not hit any black
particle.
\endi
As $\ve\to 0$, since we are putting less and less weight on energy at later times,
the single white particle that sticks to its target 
is $X_{N(\ve)}$, with $N(\ve)\to\infty$.  
In the limit, we obtain a weak solution where all black particles stick 
to each other, but all white particles hit the black particles
without sticking.    More precisely, 
for $t\in [t_k, \, t_{k-1}[$ this limit solution contains  $k-1$ 
black particles with masses $m_1,\ldots, m_{k-1}$, located
at 
$$x_j(t) =(\bar x_j+ t\bar v_j) \bfe_1\,,\qquad\qquad j=1,\ldots, k-1,$$ and one 
lumped black particle with 
mass $m_k^*= \sum_{j\ge k} m_j$, located at $b_k(t)\bfe_1$.
In addition, it contains countably many white particles with masses $M_i$, located at 
$$X_i(t)~=~\ov X_i+ t\ov V_i~=~b_i(\tau_i)\bfe_1 - (t-\tau_i)\bfe_2\,.$$
This is not a sticky solution. 
\v
{\bf Remark 3.} Given a weak solution consisting of countably many particles
$Y_i(t)$, on can introduce a measure of ``non-stickiness" 
 by setting
 $$\Phi~\doteq~\sum_{(i,j)\in NS} m_i m_j\,.$$ 
Here $m_i$ denotes the  mass of the $i$-th particle, while 
 $$NS~\doteq~\{ (i,j): \quad Y_i(t_0)= Y_j(t_0) ~\hbox{but}~
 Y_i(t)\not= Y_j(t)~~\hbox{for some times}~t>t_0\}$$
describes all couples of particles that hit each other
without sticking.

As $\ve\to 0+$, for the sequence of weak solutions considered in Remark~2
 the measure $\Phi_\ve$ of non-stickiness approaches zero.  However,  
this does not imply that the limit solution should be sticky.

\section{Non-existence and non-uniqueness for continuous initial data}
In this last section we extend Example~3 and Example~4 to the case of continuous initial data. Our main goal  is to prove:
\v
{\bf Theorem 1.} {\it  In any dimension $n\geq 2$ 
there exists an initial datum $(\bar \rho,\bar  v)$,
such that the system \eqref{SP} does not admit any sticky weak  solution
in the sense of Definition~3, not even locally in time. 
Here  $\bar \rho\in \C_c^\infty(\R^n)$ 
is the density of a probability measure, while 
$\bar v\in \C_c(\R^n)$ is a  continuous initial velocity.}
\v
{\bf Proof.} The main idea is to modify the initial data in Example~4,
 replacing each point mass, say located at $\bar y_j$,
by a smooth distribution of mass supported 
on a small ball $B_j=\{x:\, |x-\bar y_j|\leq r_j\}$. By choosing 
an appropriate initial velocity, after a very short 
time all the mass initially contained inside $B_j$  collapses to a point mass and then continues its motion as in the previous example. 
However,
a difficulty arises because the black particles move horizontally along the $x_1$ axis, while the white particles have  velocity $-\bfe_2$.  This would determine 
a discontinuity in the velocity field at the origin.
To avoid this, we need to add a horizontal component to the 
velocities of the white particles $X_k$, as shown in Fig.~\ref{f:z107}. 
We now describe the construction in greater detail.
\v
{\bf 1.}
In this step we replace the black particles $x_k$ with a smooth distribution of mass.
Choose $0<\beta<\gamma<1$ and consider the initial points and velocities
\bel{xv}\bar x_k~=~\beta^k\bfe_1\,,\qquad\qquad \bar v_k ~=~(1-\gamma^k)\bfe_1\eeq
as in (\ref{id6}).
For each  $k\geq 1$, we choose $a_k, r_k>0$  and define the smooth function 
$$\psi_k(x)~=~\left\{\bega{cl} \ds a_k
\exp\left\{{-1\over r_k^2 -|x-\bar x_k|^2}\right\}
\qquad &\hbox{if}~~|x-\bar x_k|<r_k\,,\cr\cr
0\qquad &\hbox{otherwise.}\enda\right.$$
Notice that $\psi_k\in \C^\infty_c(\R^n)$, with support contained in 
the ball $B_k\doteq\{x\in \R^n\,;~ |x-\bar x_k|\leq r_k\}$.

As initial velocity we choose a continuous function $\bar v$ such that
\bel{brdef}\bar v(x)~=~ (1-\gamma^k) \bfe_1 + {\bar x_k-x\over \sqrt{r_k}}
\qquad \quad \hbox{if}~~|x-\bar x_k|\leq r_k\,,~~~ k\geq 1.\eeq
Notice that,  for $t\geq \sqrt{r_k}$
all the mass initially located inside the ball $B_k$ gets concentrated at the single point
$$x_k(t)~=~ \bar x_k + t\bar v_k~=~[\beta^k+ t(1-\gamma^k)]\bfe_1.$$
The choice of the coefficients $a_k, r_k$ is made in two stages.

First we choose the sequence of radii $r_k\downarrow 0$ decreasing to zero fast enough
so that (i) the balls $B_k$, $k\geq 1$ are mutually disjoint and (ii) 
during the time interval $[0, \sqrt{r_k}]$, particles originating from the 
ball $B_k$ do not interact with any other particles from different balls.

Afterwards, we choose the coefficients $a_k\downarrow 0$ decreasing to zero fast enough
so that the function $\rho(x)\doteq\sum_{k\geq 1}\psi_k(x)$ is in $\C^\infty_c$.
At this stage we also observe that the conclusion of Lemma~2 remains valid
if the particle masses, instead of $m_k=\alpha^k$, are given by
$m_k~=~\int\psi_k(x)\, dx$.   Indeed, the only relevant assumption 
is that $m_k/m_{k-1}\to 0$ fast enough.
\v
{\bf 2.}
In this step we replace the white particles $X_k$ with a smooth distribution of mass.
For this purpose, it is worth noting that in Example~4
there is a lot of freedom in the choice 
of the positions and masses of the particles $X_k$.   Indeed, the only 
thing that matters is the identity $X_k(\tau_k)= b_k(\tau_k)\bfe_1$.

Let the functions $b_k(\cdot)$ and the times $\tau_k$ be as in (\ref{bar'}).
We can then consider a sequence of particles $X_k$, $k\geq 1$, 
with initial  velocity  and position given respectively by 
\bel{XV}\ov V_k~=~\bfe_1 - {\bfe_2\over k}\,,\qquad\qquad
\ov X_k~=~b_k(\tau_k)\bfe_1 - \tau_k \ov V_k\,.\eeq
Observe that (\ref{xv}) and (\ref{XV}) imply
\bel{lim0}
\lim_{k\to\infty} \bar x_k~=~\lim_{k\to\infty} \ov X_k~=~0,\qquad\qquad
\lim_{k\to\infty} \bar v_k~=~\lim_{k\to\infty} \ov V_k~=~\bfe_1\,.\eeq
We can now replace the countably many particles $X_k$ with a continuous
distribution of mass, as in the previous step.
For each $k\geq 1$, choose $\tilde a_k , R_k>0$ and define the smooth function 
$$\tilde\psi_k(x)~=~\left\{\bega{cl} \ds \tilde a_k 
\exp\left\{{-1\over R_k^2 -|x-\ov X_k|^2}\right\}
\qquad &\hbox{if}~~|x-\ov X_k|\le R_k\,,\cr\cr
0\qquad &\hbox{otherwise.}\enda\right.$$
Notice that $\tilde\psi_k\in \C^\infty_c(\R^n)$, with support contained inside the ball
$\Tilde B_k\doteq\{x\in \R^n\,;~ |x-\ov X_k|\leq R_k\}$.

As initial velocity we choose a continuous function $\bar v:\R^n\mapsto \R^n$, with
bounded support,
that satisfies (\ref{brdef})
together with 
\bel{bv3}\bar v(x)~=~ (1-\gamma^k) \bfe_1 + {\ov X_k-x\over  \sqrt{R_k}}\qquad\quad 
\hbox{if}~~|x-\ov X_k|\leq R_k\,,~~k\geq 1\,.\eeq
Notice that,  for $t\geq \sqrt{R_k}$
all the mass initially located inside the ball $\Tilde B_k$ 
gets concentrated at the single point
$$x_k(t)~=~\ov X_k + t \ov V_k \,.$$
The choice of the coefficients $\tilde a_k , R_k$ is made in two stages.

First we choose the sequence of radii $R_k\downarrow 0$ decreasing to zero fast enough
so that (i) all  the balls $B_j$, $\Tilde B_k$,  $j,k\geq 1$ are mutually disjoint and (ii) 
during the time interval $[0, \sqrt{R_k}]$, particles originating from the 
ball $\Tilde B_k$ do not interact with any other particles from different balls.

Afterwards, we choose the coefficients $\tilde a_k \downarrow 0$ 
decreasing to zero fast enough
so that the function $\tilde\rho(x)\doteq\sum_{k\geq 1}\tilde\psi_k(x)$ is in $\C^\infty_c$.
\v
{\bf 3.} 
By the same argument introduced in Example~4, for the initial data consisting 
of countably many point masses $\bar x_k, \ov X_k$ with initial velocities 
$\bar v_k, \ov V_k$, no sticky solution exists, not even locally in time.   
By the above construction, the same conclusion 
holds for an initial distribution of mass with smooth density
$$\bar \rho(x)~\doteq~C\cdot \left(\sum_{k=1}^\infty \psi_k(x) + \sum_{k=1}^\infty \tilde\psi_k(x)
\right)$$
and continuous initial velocity field
$ \bar v(\cdot)$
Here $C>0$ is a normalizing constant, chosen so that $\int\bar \rho(x)\, dx =1$.
\endproof

\v
In an entirely similar way, one can modify the initial data in   Example~3 
and obtain

{\bf Theorem 2.} {\it  In any dimension $n\geq 2$
there exists an initial datum $(\bar \rho, \bar v)$
such that the  system \eqref{SP} admits 
two distinct sticky weak  solutions. Here  $\bar \rho\in \C_c^\infty(\R^n)$ 
is the density of a probability measure, while 
$\bar v\in \C_c(\R^n)$ is a continuous initial velocity.}

\begin{figure}
\centering
\includegraphics[scale=0.6]{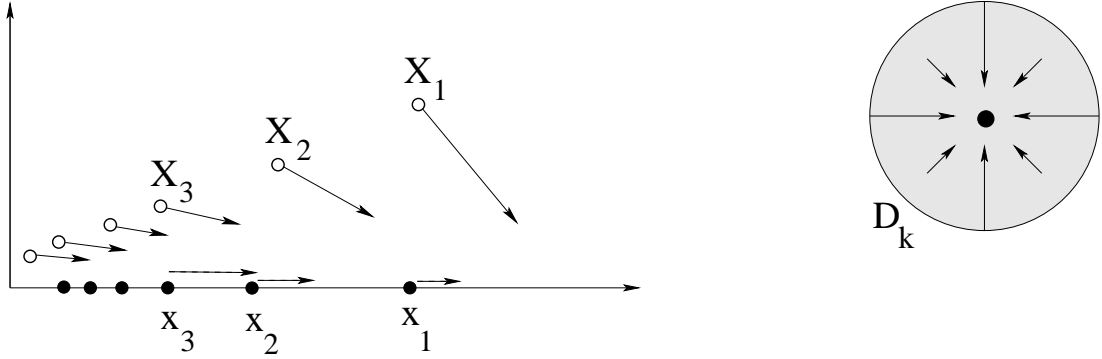}
\caption{\small Modifying the initial data in Figure~\ref{f:z106} 
in order to obtain a continuous distribution of initial velocities.
Left: the initial speed of the particles $x_k$, $X_k$ approaches the same limit as 
$k\to\infty$ and $\bar x_k, \ov X_k\to 0$.   Right:  a point mass is replaced by 
a continuous distribution on a ball $B_k$, choosing the initial velocity so that
after a short time all the mass is concentrated at one single point.
}
\label{f:z107}
\end{figure}

{\bf Remark 4.} With a more careful construction,
 counterexamples to the existence and uniqueness could be achieved
with an initial velocity distribution $\bar v\in \C_c^{1-\ve}(\R^n)$ 
which is  H\"older continuous, with any exponent strictly smaller than 1.
However, this initial velocity field cannot be 
Lipschitz continuous: in order that all the mass initially inside $B_k$ or $\Tilde B_k$ collapse 
to a point within time $t_k\to 0$, the Lipschitz constant of $\bar v$ 
in (\ref{brdef}) and
(\ref{bv3})
must tend to infinity as $k\to\infty$.   

On the other hand, if the initial velocity 
field $\bar v$ is continuous with Lipschitz constant $L$, then it is easy to see that
the Cauchy problem (\ref{SP})-(\ref{ic}) has a unique local 
solution defined for $0\leq t<L^{-1}$.
Indeed,
for each $0\leq t< L^{-1}$
the identity
$$v(t,\, x+t \bar v(x)) ~=~\bar v(x)$$
uniquely defines a Lipschitz continuous vector field $v(t,\cdot)$. Inserting
this function $v(t,\cdot)$ in (\ref{SP}), 
one obtains a linear transport equation for $\rho$, with Lipschitz velocity field and
hence with a unique solution.
\v

 \end{document}